\documentclass[12pt]{article}
\usepackage{amssymb,amsfonts,amsmath}
\usepackage{amsthm}
\usepackage{color}
\usepackage[colorlinks=true,citecolor=black,linkcolor=black,urlcolor=blue]{hyperref}

\newcommand{\arxiv}[1]{\href{http://arxiv.org/abs/#1}{\texttt{arXiv:#1}}}

\def\pf{\noindent {\it Proof.} }
\def\qed{\hfill $\square$}

\def\L{\mathcal{L}}
\def\N{\mathcal{N}}

\makeatletter\@addtoreset{figure}{section}\makeatother

\newcommand{\captionfonts}{\small}
\makeatletter
\long\def\@makecaption#1#2{
   \vskip 10\p@
   \setbox\@tempboxa\hbox{\captionfonts {#1} \ \ #2}
   \ifdim \wd\@tempboxa >\hsize
       \captionfonts {#1}\ \ #2\par
   \else
       \hbox to\hsize{\hfil\box\@tempboxa\hfil}
   \fi}
\makeatother

\newtheorem{prop}{Proposition}[section]

\newtheorem{lemma}{Lemma}[section]

\makeatletter \@addtoreset{equation}{section}
\@addtoreset{theo}{section} \makeatother

\allowdisplaybreaks

\title{A determinant representation for generalized ballot and Fuss-Catalan numbers}
\author{James J.Y. Zhao\\
\small Dongling School of Economics and Management,\\
\small University of Science and Technology Beijing,\\
\small Beijing 100083, P.R. China\\
\small and\\
\small Center for Applied Mathematics,\\
\small Tianjin University, \\
\small Tianjin 300072, P.R. China\\
\small \texttt{Email:\,zhaojy@ustb.edu.cn}}

\date{\small{December 11, 2013}}

\begin{document}

\maketitle

\begin{abstract}
In this note we introduce a determinant and then give its evaluating formula. The determinant turns out to be a generalization of the well-known ballot and Fuss-Catalan numbers, which is believed to be new. The evaluating formula is proved by showing that the determinant coincides with the number of lattice paths with $(1,0), (0,1)$-steps in the plane that stay below a boundary line of rational slope.
\end{abstract}


\section{Introduction}\label{sec-1}
For any given integers $m, n, u$ and $k$ satisfying $u\geq 0$, $k\geq 2$, $n\geq 2$, and $\ m\geq\mathrm{max}\{u+1,(k-1)(n-1)\}$, let $D(m,n,u,k)$ be a determinant as follows:
\begin{equation}\label{def-1}
D(m,n,u,k):=\det(t_{ij})_{i,j=1,2,\ldots,n-1},
\end{equation}
where
$$
t_{ij}=\binom{m-\mathrm{max}\{u,(k-1)j\}}{1-i+j}.
$$
As we shall soon see, the determinant $D(m,n,u,k)$ is a generalization of the renowned ballot and Fuss-Catalan numbers by the following proposition.
\begin{prop}\label{Dmnuk} For any given integers $m, n, u$ and $k$ satisfying $u\geq 0$, $k\geq 2$, $n\geq 2$ and $m\geq\mathrm{max}\{u+1,(k-1)(n-1)\}$ we have
\begin{align}\label{Dmnua}
D(m,n,u,k)
=\sum\limits_{i=0}^{\lfloor\frac{u}{k}\rfloor}(-1)^i\,\frac{m-(k-1)(n-1)}{m+n-1-ki}
\binom{m+n-1-k i}{n-1-i}\binom{u-(k-1)i}{i},
\end{align}
where $\lfloor x\rfloor$ is the largest integer not greater than $x$.
\end{prop}

The formula \eqref{Dmnua}, the main result of this paper, is obtained by showing that the determinant $D(m,n,u,k)$ satisfies the same recurrence and initial conditions as $|\L(u+1,1;m,n;k)|$, the cardinality of $\L(u+1,1;m,n;k)$, where $\L(u+1,1;m,n;k)$ is the set of lattice paths from $(u+1,1)$ to $(m,n)$ with $(1,0),(0,1)$-steps in the plane that stay below the line $y=\frac{x-1}{k-1}+1$.
It should be mentioned that $m$ and $n$ are indeed the variants concerning the recurrence and initial conditions, yet $u$ and $k$ are just parameters.

When $u=0$, by \eqref{Dmnua} the determinant $D(m+1,n+1,0,k+1)=\frac{m+1-kn}{m+1}\binom{m+n}{n}$ is just the solution of the generalized ballot problem which had been studied by Barbier \cite{Barbier} (cf. \cite{Gou-Ser}). Also, $D(m+1,n+1,0,2)=\frac{m-n+1}{m+1}\binom{m+n}{n}$, called the ballot number, is the solution of the original ballot problem first introduced by Bertrand \cite{Bertrand}.

Furthermore, by \eqref{Dmnua} there follows $D((k-1)n+1,n+1,0,k)=\frac{1}{(k-1)n+1}\binom{kn}{n}$, which is known as the Fuss-Catalan number \cite{Fuss91} (cf. \cite{BacKra11} or \cite[Eq. (7.67)]{G-K-P}). It also enumerates $(k-1)$-Dyck paths and $k$-ary trees. Specifically, when $k=2$ the Fuss-Catalan number reduces to $D(n+1,n+1,0,2)=\frac{1}{n+1}\binom{2n}{n}$, which is the well known Catalan number \cite{Stanley1999}, and when $k=3$ it becomes $D(2n+1,n+1,0,3)=\frac{1}{2n+1}\binom{3n}{n}$, which counts $2$-Dyck paths and ternary trees \cite{Aval08}.

The ballot number, as well as the Fuss-Catalan number, generalizes the Catalan number. These fascinating numbers have many interesting interpretations in algebra, combinatorics, and probability.
One of the most famous examples is, the classical Chung-Feller Theorem \cite{Chung-Feller} offers a graceful perspective for counting the Catalan number.
More combinatorial and algebraic interpretations of Catalan number were shown by Stanley \cite{Stanley1999, Stanley2013}. Various applications of Catalan number were collected by Koshy \cite{Koshy}.
For recent inspiring works involving generalizations of the Catalan number and the ballot number, with many references, see Krattenthaler \cite{Krat}, Sagan and Savage \cite{SagSav}, Gorsky and Mazin \cite{GorMaz}, He \cite{He}, and Bousquet-M\'{e}lou, Chapuy, and Pr\'{e}ville-Ratelle \cite{BCPR}.

Although numerous generalized ballot and Catalan numbers had been given, the determinant representation \eqref{def-1}, to the best of our knowledge, is new.
The parameter $u$ and $k$ indeed make sense while calculating certain ruin probability in the framework of classical compound binomial risk model, where the initial capital of an insurer is denoted by $u$ and each claim amount is assumed to be $k$. This shall be discussed in a following paper \cite{Zhaoac}. In particular, when $k=2$, it reduces to the classical gambler's ruin problem.

In this note, we first show a formula for $|\L(u+1,1;m,n;k)|$ in Section \ref{sec-lp}. In Section \ref{sec-pro} we give the recurrences and initial conditions for $|\L(u+1,1;m,n;k)|$ and $D(m,n,u,k)$. Then we
complete the proof of Proposition \ref{Dmnuk} by showing that $D(m,n,u,k)$ and $|\L(u+1,1;m,n;k)|$ coincide.

\section{Enumeration formula for $|\L(u+1,1;m,n;k)|$}\label{sec-lp}

In this section we shall prove the formula for $|\L(u+1,1;m,n;k)|$.

\begin{lemma}\label{lat-u1} For any given integers $m, n, u$ and $k$ satisfying $u\geq 0$, $k\geq 2$, $n\geq 2$, and $m\geq \max\{u+1,(k-1)(n-1)\}$, we have
{\small
\begin{align}\label{enu-u1}
|\L(u+1,1;m,n;k)|=
\sum\limits_{i=0}^{\lfloor\frac{u}{k}\rfloor}(-1)^i\,\frac{m-(k-1)(n-1)}{m+n-1-k i}
\binom{m+n-1-ki}{n-1-i}\binom{u-(k-1)i}{i}.
\end{align}
}
\end{lemma}

\pf
Given integers $u\geq 0$, $k\geq 2$, and $n\geq 2$, without loss of generality define $|\L(u+1,1;(k-1)(n-1),n;k)|=0$ since the point $((k-1)(n-1),n)$ stays above the boundary line $y=\frac{x-1}{k-1}+1$. Clearly \eqref{enu-u1} works when $m=(k-1)(n-1)$.

It remains to show that \eqref{enu-u1} holds when $m\geq \max\{u,(k-1)(n-1)\}+1$.
Here we need a known result on lattice paths enumeration.
In \cite[Theorem 2.1]{Zhao}, the number of lattice paths from $(a,b)$ to $(m,n)$ that stay above the line $y=kx$ was shown to be
\begin{align}\label{enu-labmnkc}
\sum\limits_{i=0}^{\lfloor\frac{b-ka}{k+1}\rfloor}(-1)^i
  \frac{n+1-km}{n+1-k(a+i)}\binom{m+n-(k+1)(a+i)}{m-a-i}\binom{b-k(a+i)}{i},
\end{align}
where $a, b, m, n$ and $k$ are integers satisfying $k\geq 1$, $0\leq a\leq m$, $ka\leq b\leq n$ and $n\geq km$, by mathematical induction based on a forward recursion.

To obtain \eqref{enu-u1} from \eqref{enu-labmnkc}, first reflect each path counted by \eqref{enu-labmnkc} along the diagonal $y=x$ in the plane. The resulting path obviously belongs to the set of lattice paths from $(b,a)$ to $(n,m)$ that stay below the line $y=x/k$ which is denoted by $\N_{1/k}(b,a;n,m)$.
Make the substitutions $a\rightarrow b, b\rightarrow a, m\rightarrow n$, and $n\rightarrow m$ in \eqref{enu-labmnkc}. Then we have
\begin{align}\label{enu-labmna}
|\N_{1/k}(a,b;m,n)|\hskip-2pt=\hskip-3pt\sum\limits_{i=0}^{\lfloor\frac{a-kb}{k+1}\rfloor}(-1)^i
  \frac{m+1-kn}{m+1-k(b+i)}\binom{m+n-(k+1)(b+i)}{n-b-i}\binom{a-k(b+i)}{i},
\end{align}
for integers $a, b, m, n$ and $k$ satisfying $k\geq 1$, $0\leq b\leq n$, $kb\leq a\leq m$, $n\geq 1$ and $m\geq kn$.

Next for any lattice path $\xi\in\L(u+1,1;m,n;k)$, shift the $x$-axis upward for one unit and the $y$-axis rightward for one unit. The resulting path clearly belongs to the set $\N_{1/(k-1)}(u,0;m-1,n-1)$. This map is obviously reversible and unique. Thus
\begin{align*}
|\L(u+1,1;m,n;k)|=|\N_{1/(k-1)}(u,0;m-1,n-1)|.
\end{align*}
After making the substitutions $k\rightarrow k-1, a\rightarrow u, b\rightarrow 0, m\rightarrow m-1$, and $n\rightarrow n-1$ in \eqref{enu-labmna}, we obtain \eqref{enu-u1} for $m\geq \max\{u,(k-1)(n-1)\}+1$. This completes the proof.
\qed

\section{Proof of Proposition \ref{Dmnuk}}\label{sec-pro}
The aim of this section is to complete the proof of Proposition \ref{Dmnuk}. This is achieved by showing that both $D(m,n,u,k)$ and $|\L(u+1,1;m,n;k)|$ have the same recurrence and initial conditions. First we have the following result for $|\L(u+1,1;m,n;k)|$.
\begin{lemma}\label{lem-1}
Given integers $u\geq 0$ and $k\geq 2$, the numbers $|\L(u+1,1;m,n;k)|$ satisfy the recurrence
\begin{eqnarray}\label{rec-L}
\begin{array}{r}
|\L(u+1,1;m,n;k)|=|\L(u+1,1;m-1,n;k)|+|\L(u+1,1;m,n-1;k)|,\\
n\geq 3,\quad m\geq\mathrm{max}\{u+2,\,(k-1)(n-1)+1\},
\end{array}
\end{eqnarray}
with the initial conditions
\begin{eqnarray}
&\hskip -1cm |\L(u+1,1;m,2;k)|=m-\max\{u,k-1\},\ m\geq\mathrm{max}\{u+2,(k-1)(n-1)+1\},\label{eq-Lc1}\\
&\hskip -1cm|\L(u+1,1;u+1,n;k)|=1, \hskip 3.9cm  2\leq n\leq \lfloor\frac{u}{k-1}\rfloor+1,\, u\geq k-1,\label{eq-Lc2}\\
&\hskip -1cm|\L(u+1,1;(k-1)(n-1),n;k)|=0, \hskip 3.53cm n\geq\mathrm{max}\{2,\lceil\frac{u}{k-1}\rceil+1\},\label{eq-Lc3}
\end{eqnarray}
where $\lceil x\rceil$ is the least integer not less than $x$.
\end{lemma}

\pf
The recurrence \eqref{rec-L} is obviously true from the definition of $|\L(u+1,1;m,n;k)|$.

Note that the initial conditions are distributed on the lines $y=2$, $x=u+1$ (if $u\geq k-1$) and $y=\frac{x}{k-1}+1$, corresponding to \eqref{eq-Lc1}, \eqref{eq-Lc2}, and \eqref{eq-Lc3}, respectively, in the plane.

\begin{figure}[!ht]
\begin{center}
\setlength\unitlength{3.pt}
\begin{picture}(96,57)
\thinlines
\put(2,12){\vector(1,0){92}}
\put(92,9.7){\small $x$}
\put(14,4){\vector(0,1){52}}
\put(11.5,54.5){\small $y$}
\put(8,12){\line(0,1){1}}
\multiput(20,12)(6,0){11}{\line(0,1){1}}
\put(14,15){\line(1,0){1}}
\multiput(14,18)(0,6){5}{\line(1,0){1}}
\multiput(20,18)(6,0){11}{\circle*{1.0}}
\multiput(26,24)(6,0){10}{\circle*{1.0}}
\multiput(38,30)(6,0){8}{\circle*{1.0}}
\multiput(50,36)(6,0){6}{\circle*{1.0}}
\multiput(62,42)(6,0){4}{\circle*{1.0}}
\multiput(74,48)(6,0){2}{\circle*{1.0}}
\put(5.9,9.1){\footnotesize $-1$}
\put(12,9.1){\footnotesize $0$}
\put(19.2,9){\footnotesize $1$}
\put(25.2,9){\footnotesize $2$}
\put(31.2,9){\footnotesize $3$}
\put(37.2,9){\footnotesize $4$}
\put(43.2,9){\footnotesize $5$}
\put(49.2,9){\footnotesize $6$}
\put(55.2,9){\footnotesize $7$}
\put(61.2,9){\footnotesize $8$}
\put(67.2,9){\footnotesize $9$}
\put(72.5,9){\footnotesize $10$}
\put(78.5,9){\footnotesize $11$}
\put(10,14){\scriptsize $0.5$}
\put(11,17.5){\footnotesize $1$}
\put(11,23){\footnotesize $2$}
\put(11,29){\footnotesize $3$}
\put(11,35){\footnotesize $4$}
\put(11,41){\footnotesize $5$}
\put(2,9){\line(2,1){91}}
\put(87.5,50){\small $y=\frac{x}{2}+\frac{1}{2}$}
\put(25,15){\scriptsize $A$}
\put(31,15){\scriptsize $B$}
\multiput(24.5,19.3)(6,0){10}{\scriptsize $1$}
\put(24.5,25.3){\scriptsize $0$}
\put(30.5,25.3){\scriptsize $1$}
\put(36.5,25.3){\scriptsize $2$}
\put(42.5,25.3){\scriptsize $3$}
\put(48.5,25.3){\scriptsize $4$}
\put(54.5,25.3){\scriptsize $5$}
\put(60.5,25.3){\scriptsize $6$}
\put(66.5,25.3){\scriptsize $7$}
\put(72.5,25.3){\scriptsize $8$}
\put(78.5,25.3){\scriptsize $9$}
\put(36.5,31.3){\scriptsize $0$}
\put(42.5,31.3){\scriptsize $3$}
\put(48.5,31.3){\scriptsize $7$}
\put(53.5,31.3){\scriptsize $12$}
\put(60,31.3){\scriptsize $18$}
\put(66,31.3){\scriptsize $25$}
\put(72,31.3){\scriptsize $33$}
\put(78,31.3){\scriptsize $42$}
\put(48,37.3){\scriptsize $0$}
\put(54,37.3){\scriptsize $12$}
\put(60,37.3){\scriptsize $30$}
\put(65.5,37.3){\scriptsize $55$}
\put(72,37.3){\scriptsize $88$}
\put(77.5,37.3){\scriptsize $130$}
\put(59.5,43.3){\scriptsize $0$}
\put(65.5,43.3){\scriptsize $55$}
\put(71.5,43.3){\scriptsize $143$}
\put(77.5,43.3){\scriptsize $273$}
\put(71.5,49.3){\scriptsize $0$}
\put(77.5,49.3){\scriptsize $273$}
\thicklines
\multiput(32,18)(6,6){2}{\line(1,0){6}}
\multiput(38,18)(6,6){2}{\line(0,1){6}}
\multiput(44,30)(12,6){3}{\line(1,0){6}}
\multiput(50,30)(12,6){3}{\line(1,0){6}}
\multiput(56,30)(12,6){2}{\line(0,1){6}}
{\color{blue}
\linethickness{0.5mm}
\put(26,18){\line(1,0){6}}
}
\end{picture}
\end{center}
\vskip -1cm
\caption{A lattice path $\xi\in\L(2,1;11,5;3)$ where $1=u<k-1=2$. In this case the boundary line is $y=\frac{x}{2}+\frac{1}{2}$. The blue portion $AB$ is the first $k-1-u$ step, which has to be horizontal. The number above each point $(m,n)$ is equal to $|\L(2,1;m,n;3)|$.\label{fig-1}}
\end{figure}
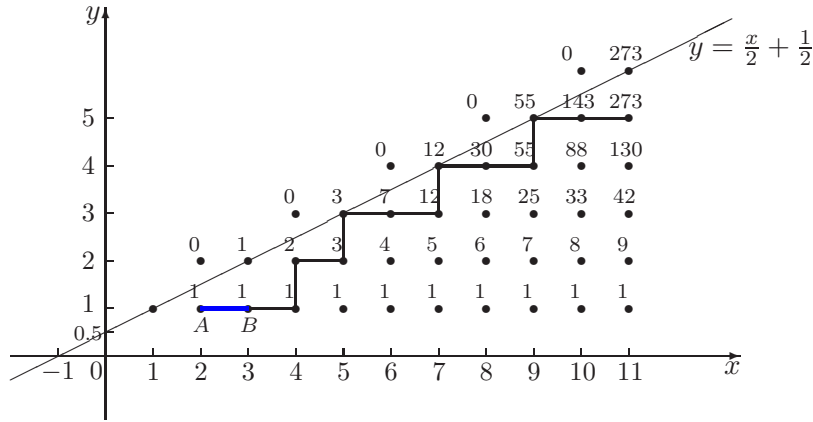

When $n=2$, if $0\leq u< k-1$, for each lattice path $\xi\in \L(u+1,1;m,2;k)$ with $m>u+1$, the first $k-1-u$ steps must be horizontal because of the restriction that $\xi$ must stay below the line $y=\frac{x-1}{k-1}+1$. For example see the blue portion $AB$ in Figure \ref{fig-1}. So $|\L(u+1,1;m,2;k)|=|\L(k,1;m,2;k)|=\binom{m-k+2-1}{2-1}=m-(k-1)$. If $u\geq k-1\geq 1$, by the definition of $\L(u+1,1;m,n;k)$ it is easy to see that $|\L(u+1,1;m,2;k)|=\binom{m-(u+1)+2-1}{2-1}=m-u$. Hence \eqref{eq-Lc1} follows.

The initial conditions \eqref{eq-Lc2} and \eqref{eq-Lc3} are easy to obtain by the definition of $\L(u+1,1;m,n;k)$. The restrictions $2\leq n\leq \lfloor\frac{u}{k-1}\rfloor+1$ and $u\geq k-1$ in \eqref{eq-Lc2} are required to make sure that each lattice path stays below the boundary line $y=\frac{x-1}{k-1}+1$. Note that when $u<k-1$, the initial conditions consist of only \eqref{eq-Lc1} and \eqref{eq-Lc3}. See for example in Figure \ref{fig-1}.

Now we have shown that \eqref{rec-L}--\eqref{eq-Lc3} hold. This completes the proof.
\qed

Next we shall prove the recurrence and initial conditions for $D(m,n,u,k)$.

\begin{lemma}\label{lem-2}
Given integers $u\geq 0$ and $k\geq 2$, the determinant $D(m,n,u,k)$ satisfies the recurrence
\begin{eqnarray}\label{eq-reD}
\begin{array}{r}
D(m,n,u,k)=D(m-1,n,u,k)+D(m,n-1,u,k),\\
n\geq 3,\ m\geq\mathrm{max}\{u+2,(k-1)(n-1)+1\},
\end{array}
\end{eqnarray}
with the initial conditions
\begin{eqnarray}
&\hskip -0.5cm D(m,2,u,k)=m-\mathrm{max}\{u,k-1\},\hskip 0.33cm m\geq\mathrm{max}\{u+2,(k-1)(n-1)+1\},\label{eq-m2}\\
&\hskip -0.5cm D(u+1,n,u,k)=1,\hskip 4.2cm 2\leq n\leq \lfloor\frac{u}{k-1}\rfloor+1,\, u\geq k-1,\label{eq-u1}\\
&\hskip -0.5cm D((k-1)(n-1),n,u,k)=0,\hskip 3.8cm n\geq\mathrm{max}\{2,\lceil\frac{u}{k-1}\rceil+1\}.\label{eq-du0}
\end{eqnarray}
\end{lemma}

\pf
To prove \eqref{eq-reD}, we first need to perform an elementary row operation on the matrix $(t_{ij})_{i,j=1,2,\ldots,n-1}$ defined in \eqref{def-1}. That is, replace the $i$-th row by adding the $(i+1)$-th row to it in increasing order of $i$. Let $(t^\ast_{ij})_{i,j=1,2,\ldots,n-1}$ denote the resulting matrix. Clearly $D(m,n,u,k)=\det(t^\ast_{ij})_{i,j=1,2,\ldots,n-1}$, where
$$
t^\ast_{ij}=\left\{
\begin{array}{lr}
\binom{m+1-\mathrm{max}\{u,(k-1)j\}}{1-i+j},&\quad i=1,2,\ldots,n-2,j=1,2,\ldots,n-1,\\[7pt]
\binom{m-\mathrm{max}\{u,(k-1)j\}}{1-i+j},&\quad i=n-1,j=1,2,\ldots,n-1.
\end{array}
\right.
$$

The first $n-2$ rows of the matrix $(t^\ast_{ij})_{i,j=1,2,\ldots,n-1}$ are obtained easily by applying the well known identity $\binom{n}{k}+\binom{n}{k-1}=\binom{n+1}{k}$ for integer $k$ \cite[Eq. (5.8)]{G-K-P}, and it is clear that both $(t_{ij})_{i,j=1,2,\ldots,n-1}$ and $(t^\ast_{ij})_{i,j=1,2,\ldots,n-1}$ have the same bottom row, $\left(0,\ldots,0,1,m-\mathrm{max}\{u,(k-1)(n-1)\}\right)$, whose zero entries are obtained by the definition $\binom{n}{k}=0$ when $n\geq 0$ and $k<0$.

Now let $t_j^\ast$ denote the $j$-th column of the matrix $(t_{ij}^\ast)_{i,j=1,2,\ldots,n-1}$, and let $^\mathrm{T}$ denote transpose. Note that $t_{n-1}^\ast=\alpha+\beta$, where
$\alpha=(t_{1,n-1}^\ast,\ldots,t_{n-2,n-1}^\ast,m+1-{\rm max}\{u,(k-1)(n-1)\})^{\mathrm T}$ and $\beta=(\underbrace{0,\ldots,0}_{(n-2)'s},-1)^{\rm T}$.
Thus by basic properties of determinant and the definition of $D(m,n,u,k)$, we have
\begin{align}\label{eq-A-B}
D(m,n,u,k)&=\det\left(t^\ast_1,\ldots,t^\ast_{n-2},\alpha+\beta\right)\nonumber\\
 &=\det\left(t^\ast_1,\ldots,t^\ast_{n-2},\alpha\right)
  +\det\left(t^\ast_1,\ldots,t^\ast_{n-2},\beta\right)\nonumber\\
 &=D(m+1,n,u,k)-D(m+1,n-1,u,k),
\end{align}
for $n\geq 3$ and $m\geq \mathrm{max}\{u+1,(k-1)(n-1)\}$.
Letting $m\rightarrow m-1$ in \eqref{eq-A-B} gives \eqref{eq-reD}.

Next we shall prove the initial conditions of $D(m,n,u,k)$.
Given integers $u\geq 0$ and $k\geq 2$, first by the definition of $D(m,n,u,k)$ we have $D(m,2,u,k)=\left|\binom{m-\mathrm{max}\{u,k-1\}}{1}\right|=m-\mathrm{max}\{u,k-1\}$ when $m\geq\mathrm{max}\{u+2,(k-1)(n-1)+1\}$, which leads to \eqref{eq-m2}.

When $m=u+1$, $2\leq n\leq \lfloor\frac{u}{k-1}\rfloor+1$, and $u\geq k-1$, it follows that $t_{ij}=\binom{u+1-\mathrm{max}\{u,(k-1)j\}}{1-i+j}=\binom{1}{1-i+j}$ since $u\geq (k-1)(n-1)\geq (k-1)j$ for $1\leq j\leq n-1$. Then we have
$$
t_{ij}=\left\{
\begin{array}{lr}
0,&\quad 1\leq i<j\leq n-1,\\
1,&\quad i=j,\\
1,&\quad i=j+1,\\
0,&\quad i>j+1.
\end{array}
\right.
$$
Clearly $(t_{ij})_{i,j=1,2,\ldots,n-1}$ is a lower triangular matrix with $1$'s on the main diagonal when $m=u+1$, $2\leq n\leq \lfloor\frac{u}{k-1}\rfloor+1$, and $u\geq k-1$, therefore \eqref{eq-u1} follows.

Note that if $u<k-1$, by the restriction $m\geq\mathrm{max}\{u+1,(k-1)(n-1)\}$, the variant $m$ should start at $(k-1)(n-1)$ for $k,n\geq 2$.

Finally when $u>0$ and $n\geq \lceil\frac{u}{k-1}\rceil+1$, it is clear that $(k-1)(n-1)\geq u$. Let us consider the matrix corresponding to $D((k-1)(n-1),n,u,k)$. By the definition of $D(m,n,u,k)$ the entry $t_{i,n-1}=\binom{(k-1)(n-1)-{\rm max}(u,(k-1)(n-1))}{1-i+n-1}=\binom{0}{n-i}=0$ for $i=1,2,\ldots,n-1$, which means that $D((k-1)(n-1),n,u,k)=0$.
When $u=0$ and $n\geq 2$, a similar discussion also leads to the same outcome. Then we have \eqref{eq-du0} as desired.

This completes the proof.
\qed

We are now in a position to complete the proof of Proposition \ref{Dmnuk}.

\noindent {\it Proof of Proposition \ref{Dmnuk}.}
Let $m,n,u$ and $k$ be integers where $u\geq 0$, $k\geq 2$, $n\geq 2$ and $m\geq\mathrm{max}\{u+1,\,(k-1)(n-1)\}$. First by Lemma \ref{lem-1} and Lemma \ref{lem-2}, it is clear that $D(m,n,u,k)$ satisfies the same recurrence and initial conditions as $|\L(u+1,1;m,n;k)|$, so they agree.
That is, $D(m,n,u,k)=|\L(u+1,1;m,n;k)|$. Together with Lemma \ref{lat-u1}, there follows \eqref{Dmnua}.
This completes the proof.
\qed

\vspace{.2cm} \noindent{\bf Acknowledgments}
I am very grateful to Arthur L.B. Yang and Guo-Ce Xin for inspiring discussions and valuable comments on an earlier version of this paper. This work was supported by the NSF of China Grant 71171018, the Fundamental Research Funds for the Central Universities Grant FRF-BR-12-008, and the Specialized Research Fund for the Doctoral Program of Higher Education No. 20130006110001.


\begin{thebibliography}{99}

\bibitem{Aval08}
J.-C. Aval, Multivariate Fuss-Catalan numbers, Discrete Math. 308(20):4660--4669, 2008.

\bibitem{BacKra11}
R. Bacher and C. Krattenthaler, Chromatic statistics for triangulations and Fu{\ss}-Catalan complexes, Electron. J. Combin., 18(1) \#P152, 2011.

\bibitem{Barbier}
E. Barbier, Note: Calcul des probabilites. Generalisation du probleme resolu par M. J. Bertrand, Co. Re. Acad. Sci. Paris, 105, p. 407, 1887.

\bibitem{Bertrand}
J. Bertrand, Solution d'un probl\'{e}me, Comptes Rendus de l'Acad\'{e}mie des Sciences, Paris, 105  p. 369. 1887.

\bibitem{BCPR}
M. Bousquet-M\'{e}lou, G. Chapuy, and L.-F. Pr\'{e}ville-Ratelle, The representation of the symmetric group on $m$-Tamari intervals, Adv. Math., 247(10):309--342, 2013.

\bibitem{Chung-Feller}
K.L. Chung, and W. Feller, On fluctuations in coin-tossing, Proc. Natl. Acad. Sci. U.S.A., 35(10):605--608, 1949.

\bibitem{He}
T.X. He, Parametric Catalan numbers and Catalan triangles, Linear Algebra Appl., 438(3):1467--1484, 2013.

\bibitem{Fuss91}
N. Fuss, Solutio qu{\ae}stionis, quot modis polygonum n laterum in polygona m laterum, per diagonales resolvi qu{\ae}at, Nova acta academi{\ae} scientiarum Imperialis Petropolitan{\ae}, 9:243--251, 1791.

\bibitem{GorMaz}
E. Gorsky and M. Mazin, Compactified Jacobians and $q,t$-Catalan numbers, I, J. Combin. Theory Ser. A, 120(1):49-63, 2013.

\bibitem{Gou-Ser}
I.P. Goulden and L.G. Serrano, Maintaining the spirit of the reflection principle when the boundary has arbitrary integer slope, J. Combin. Theory Ser. A, 104(2):317--326, 2003.

\bibitem{G-K-P}
R.L. Graham, D.E.  Knuth, and O. Patashnik, Concrete Mathematics: A Foundation for Computer Science, 2nd ed. Reading, MA: Addison-Wesley, 1994.

\bibitem{Koshy}
T. Koshy, Catalan Numbers with Applications, Oxford University Press, New York, 2009.

\bibitem{Krat}
C. Krattenthaler, Determinants of (generalised) Catalan numbers, J. Statist. Plann. Inference, 140(8):2260--2270, 2010.

\bibitem{SagSav}
B.E. Sagan, C.D. Savage, Mahonian pairs, J. Combin. Theory Ser. A, 119(3):526--545, 2012.

\bibitem{Stanley1999}
R.P. Stanley, Enumerative Combinatorics, Vol. 2. Cambridge University Press, Cambridge, UK, 1999.

\bibitem{Stanley2013}
R.P. Stanley, \href{http://www-math.mit.edu/~rstan/ec/catadd.pdf}{\texttt{Catalan addendum}}, version of 25 May 2013.

\bibitem{Zhao}
J.J.Y. Zhao, Koroljuk's formula for counting lattice paths revisited, preprint. \arxiv{1306.6015v1}

\bibitem{Zhaoac}
J.J.Y. Zhao, A combinatorial approach for calculating certain ruin probabilities, preprint.

\end{thebibliography}
\end{document}